\documentclass[11pt]{article}
\usepackage[english]{babel}
\usepackage{amsmath}
\usepackage{amsthm}
\usepackage{amssymb}
\usepackage[dvips]{color}
\newtheorem{theorem}{Theorem}[section]
\newtheorem{definition}[theorem]{Definition}

\newtheorem{remark}[theorem]{Remark}

\title{\bf A Brief History of Quaternions and of the Theory of Holomorphic
Functions of Quaternionic Variables}

\author{Amy Buchmann\\Department of Mathematics \\and Computer Sciences
\\Schmid College of Science\\Chapman University\\Orange, CA 92866 USA,
\\buchm100@mail.chapman.edu}

\date{March 2009 }
\begin{document}
\maketitle
\begin{abstract}
In this paper I will give a brief history of the discovery (Hamilton, 1843) of quaternions. I will
address the issue of why a theory of triplets (the original goal of Hamilton) could not be developed.
Finally, I will discuss briefly the history of various attempts to define holomorphic functions on
quaternionic variables.
\end{abstract}

\vspace{10mm}
\begin{center}
Advised by Daniele C. Struppa, PhD.\\
Chancellor \\
Chapman University \\
One University Drive \\ 
Orange, CA 92866\\
\end{center}

\newpage
\section{Introduction}

The discovery of the quaternions is one of the most well documented
discoveries in mathematics. In general, it is very rare that the date and location of
a major mathematical discovery are known. In the case of quaternions, however, we know
that they were discovered by the
Irish mathematician, William Rowan Hamilton on October 16$^{th}$, 1843 (we will see later how
we come to be so precise). 

The early 19$^{th}$ century was a very exciting time for Complex Analysis.
Though complex numbers had been discussed in works published in the 1500s,
the study of complex numbers was often dismissed as useless. It was not
until centuries later that the practicality of complex numbers was really
understood. In the late 1700s, Euler made significant contributions to
complex analysis, but most of the fundamental results which now form the core of complex analysis were
discovered between 1814 and 1851 by Cauchy, Riemann, among others \cite%
{Needham}. Hamilton himself, before discovering quaternions, was involved with complex numbers. In 1833, he completed the Theory
of Couplets, which was at the time regarded as a new algebraic representation of the Complex Numbers. He
represented two real numbers $a$ and $b$ as a couple $(a,b)$, and then defined the
additive operation to be $(a,b)+(c,d)=(a+c,b+d)$ and the multiplicative
operation to be $(a,b)\cdot (c,d)=(ac-bd,ad+bc)$ thus defining the algebraic
definition of complex numbers. From a modern perspective, this is nothing more than a well understood way to
represent complex numbers, but in 1833, the understanding of the algebraic
structure of what we now call the complex field was not fully grasped.

As his next natural, step, Hamilton wanted to extend the complex numbers
to a new algebraic structure with each element consisting of one real part
and two distinct imaginary parts. This would be known as the Theory of
Triplets. There are natural mathematical reasons why one would want to attempt such
a construction, but Hamilton was guided, as well, by a desire to use these triplets to
represent rotations in three-dimensional space, just like complex numbers could be used
to represent rotations in the two-dimensional plane.

Hamilton worked unsuccessfully at creating this algebra for over
10 years, and finally had a breakthrough on October 16$^{th}$, 1843 while on a
walk with his wife, Lady Hamilton. They had been walking along the Royal
Canal in Dublin when it occurred to Hamilton that his new algebra would
require three rather than two imaginary parts. In order to do this, he could create a
new algebraic structure consisting of one real part and three imaginary
parts i, j, and k. For this new structure to work, Hamilton realized that
these new imaginary elements would have to satisfy the following
conditions
\begin{equation}\label{imaginary}
i^{2}=j^{2}=k^{2}=ijk=-1.
\end{equation}
Hamilton carved these results on the
nearby Broome Bridge. Unfortunately the carvings no longer remain today \cite {Altmann}.
However, his discovery was so significant that every year on
October 16th, the Mathematics Department of the National University of
Ireland, Maynooth, holds a Hamilton Walk to Broome Bridge commemorating his
discovery. This sequence of events is documented in a famous letter that
Hamilton wrote to his son, which I attach in the appendix.

In this paper, I will first describe the skew field of quaternions, and I will
then attempt to explain why Hamilton had to abandon the theory of triplets. I will conclude with a
section that traces one of the most important developments in the study of quaternions since
Hamilton, namely the attempt to replicate, for quaternionic functions, the theory of
holomorphic functions that has so much importance in the study of the complex plane.

\bigskip

\section{The algebra of quaternions}

In this section, I set the stage for the rest of the paper and I provide the
basic algebraic
definitions of quaternions\cite{Gallian}. To better understand and appreciate the
discovery of quaternions, it is important to understand them as an algebraic
structure. The quaternions, often denoted by $\mathbb{H}$, in honor of their discoverer,
constitute a non-commutative field, also known as a skew
field, that extends the field $\mathbb{C}$ of complex numbers.

In modern day mathematics, we use
abstract algebra to describe algebraic structures as groups, rings, and
fields. A group, for example, is a non-empty set together with an
associative binary operation and which contains an identity, and an inverse
for every element. Galois first used the term group in around 1830, so this
was a new and exciting branch of mathematics in Hamilton's time. In abstract
algebra, a field is a set $\mathbb{F}$, together with two associative binary
operations, typically referred to as addition and multiplication. In order to have
a field, we require
that $\mathbb{F}$ is an Abelian group under addition (i.e. addition is commutative
and its identity is denoted by 0), that $\mathbb{F}\backslash \{0\}$ is an Abelian
group under multiplication, and that multiplication must distribute over
addition. In short, this guarantees the ability to perform addition,
subtraction, multiplication, and division on all elements in the field (except
division by zero). If $\mathbb{F}\backslash \{0\}$ is not an Abelian group under
multiplication, then we call $\mathbb{F}$ a skew field. In this case, we are still able to perform
addition, subtraction, multiplication, and division, but the multiplication
is not commutative. While the terminology for fields and groups was not fully established
until after Hamilton's discovery of quaternions, this abstract approach was
already being used at this time.

To begin with, the field of complex numbers is defined by
$$\mathbb{C}=\{a+ib | a,b\in \mathbb{R},  i^{2}=-1\}.$$
This means that every complex number can
be written in the form $a+ib$ where $a$ and $b$ are real numbers and $i$ is
an imaginary unit, i.e. $i^{2}=-1.$ We now want to construct a new
field, $\mathbb{H}$, such that $\mathbb{C}$ is a subset of $\mathbb{H}$, and such that
the new operations in $\mathbb{H}$ are compatible with the old operations in
$\mathbb{C}.$ This means we
are looking for a field extension of $\mathbb{C}$ according to the following definition.

\begin{definition} Given two fields, $\mathbb{E}$ and $\mathbb{F}$, we say that
$\mathbb{E}$ is a field extension over $\mathbb{F}$ if $\mathbb{F}$ is a subset of
$\mathbb{E}$ and the
operations of $\mathbb{F}$ are those of $\mathbb{E}$ restricted to $\mathbb{F}$ \cite{Gallian}.
\end{definition}

\begin{remark} For example, $\mathbb{C}$ is a field extension over $\mathbb{R}$. Indeed,
 $\mathbb{R}$
is a subset of $\mathbb{C}$ because every element of $\mathbb{R}$
can be written in the form $a+ib$ where $b=0$. Finally, the additive and
multiplicative operators of $\mathbb{R}$ are consistent with those of $\mathbb{C}$ when $b=0$.
In general, we say that $\mathbb{E}$ has degree n over $\mathbb{F}$ if $\mathbb{E}$ has
dimension n as a vector space over $\mathbb{F}$. Because $\{1,i\}$ forms a $2$ dimensional
vector space of $\mathbb{C}$ over $\mathbb{R}$, the field of complex numbers
is a degree $2$ field extension over the field of real numbers.
\end{remark}

In the specific case of quaternions, $\mathbb{H}$ is constructed by adding two new
elements $j$ and $k$ such that $i^{2}=j^{2}=k^{2}=ijk=-1.$ The field of
quaternions can then be written as

$$\mathbb{H}=\{q=q_{0}+q_{1}i+q_{2}j+q_{3}k | q_{t}\in \mathbb{R},
i^{2}=j^{2}=k^{2}=ijk=-1\}.$$

We can now begin to talk about the additive and multiplicative operations
that can be defined on elements of $\mathbb{H}$ and turn it into a field \cite{Kuipers}.
Let $p=p_{o}+p_{1}i+p_{2}j+p_{3}k$ and $q=q_{o}+q_{1}i+q_{2}j+q_{3}k$ be two
elements in $\mathbb{H}.$ We say that $p=q$ if and only if $p_{t}=q_{t}$ for all $t\in \{0,1,2,3\}.$ Addition is defined in the natural
sense, that is $p+q=p_{o}+q_{o}+(p_{1}+q_{1})i+(p_{2}+q_{2})j+(p_{3}+q_{3})k.$

It takes a little bit more work to define the multiplication of p and q.
Using the distributivity of multiplication over addition, the product

$$pq=(p_{o}+p_{1}i+p_{2}j+p_{3}k)(q_{o}+q_{1}i+q_{2}j+q_{3}k)$$

\noindent
becomes

$$
pq=p_{o}q_{o}+p_{o}q_{1}i+p_{o}q_{2}j+p_{o}q_{3}k+p_{1}q_{0}i+p_{1}q_{1}i^{2}+p_{1}q_{2}ij+$$
$$+p_{1}q_{3}ik+p_{2}q_{0}j+p_{2}q_{1}ji+p_{2}q_{2}j^{2}+p_{2}q_{3}jk+p_{3}q_{0}k+p_{3}q_{1}ki+p_{3}q_{2}kj+p_{3}q_{3}k^{2}
$$

Notice that we must be careful simplifying this step because multiplication
is commutative for real numbers but not for imaginary elements. Using the
basic properties of quaternions and identities (\ref{imaginary}) we can rewrite the multiplication again.

$$
pq=p_{0}q_{0}+p_{0}q_{1}i+p_{0}q_{2}j+p_{0}q_{3}k+p_{1}q_{0}i-p_{1}q_{1}+p_{1}q_{2}k-p_{1}q_{3}j+$$
$$+p_{2}q_{0}j
-p_{2}q_{1}k-p_{2}q_{2}+p_{2}q_{3}i+p_{3}q_{0}k+p_{3}q_{1}j-p_{3}q_{2}i-p_{3}q_{3}
$$

We regroup the terms according to the imaginary units, and we obtain

$$
pq=p_{0}q_{0}-(p_{1}q_{1}+p_{2}q_{2}+p_{3}q_{3})+(p_{o}q_{1}+p_{1}q_{0}+p_{2}q_{3}-p_{3}q_{2})i+$$
$$+(p_{0}q_{2}-p_{1}q_{3}+p_{2}q_{0}+p_{3}q_{1})j+(p_{0}q_{3}+p_{1}q_{2}-p_{2}q_{1}+p_{3}q_{0})k
$$

\bigskip We may also represent this product as the product of two matrices
where

$pq=r_{0}+ir_{1}+jr_{2}+kr_{3}$ such that $\left[
\begin{array}{c}
r_{0} \\
r_{1} \\
r_{2} \\
r_{3}%
\end{array}%
\right] =\left[
\begin{array}{cccc}
p_{0} & -p_{1} & -p_{2} & -p_{3} \\
p_{1} & p_{0} & -p_{3} & p_2 \\
p_{2} & p_{3} & p_{0} & -p_{1} \\
p_{3} & -p_{2} & p_{1} & p_{0}%
\end{array}%
\right] \left[
\begin{array}{c}
q_{0} \\
q_{1} \\
q_{2} \\
q_{3}%
\end{array}%
\right] $

\bigskip Notice that if we rewrite the quaternions p and q as $p=p_{0}+\mathbf{p}$
and $q=q_{0}+\mathbf{q}$ where $\mathbf{p}=p_{1}i+p_{2}j+p_{3}k$ and $%
\mathbf{q}=q_{1}i+q_{2}j+q_{3}k$, we can rewrite the product of $p$ and $q$
as $pq=p_{0}q_{0}-\mathbf{p}\cdot \mathbf{q}+p_{0}\mathbf{q}+q_{0}\mathbf{p}+%
\mathbf{p}\times \mathbf{q.}$

We are also able to define a complex conjugate of a quaternion $q$. Let $%
q=q_{0}+\mathbf{q,}$ then the complex conjugate of $q$, denoted $\overline{q}
$ is the quaternion $\overline{q}=q_{0}-\mathbf{q.}$ We can also define the
norm of a quaternion $q,$denoted $|q|$ to be $\sqrt{q\overline{q}}$. Because we are working in a skew field, it is important to note that $q\overline{q}=\overline{q}q=|q|^{2}.$

We can now show that the quaternions are indeed a field, which means that
every nonzero element has a multiplicative inverse. For a nonzero quaternion $q$, the
inverse of $q$, denoted $q^{-1}$, is the quaternion $q^{-1}=\frac{\overline{q%
}}{|q|^{2}}.$ We can easily verify this by preforming the multiplication $
qq^{-1}=\frac{q\overline{q}}{|q|^{2}}=\frac{|q|^{2}}{|q|^{2}}=1$ and similarly $%
q^{-1}q=\frac{\overline{q}q}{|q|^{2}}=\frac{|q|^{2}}{|q|^{2}}=1$

\textbf{\bigskip }

\section{Why dimension $4$?}

In this section, I will investigate Hamilton's breakthrough concerning the
necessity of three distinct imaginary parts. Why was Hamilton unable to
create his Theory of Triplets? What made him realize that he would need to
add a third imaginary element? This can be understood with a basic
understanding of field extensions.

If Hamilton had been able to develop his Theory of Triplets, he would have
effectively built a degree three field extension of $\mathbb{R}$ whose
vector space forms the basis $\{1,i,j\}$ over $\mathbb{R}$ such that $%
i^{2}=j^{2}=-1$.

Let us then call $\mathbb{D}$ the hypothetical field generated by the elements
$\{1,i,j\}$ satisfying $i^{2}=j^{2}=-1.$ Because $\mathbb{D}$ is a field, it must be closed under
multiplication. Thus the multiplication of $i$ and $j$ must result in some
element that is already in the field. In order to find a contradiction, we
investigate all of the possibilities of the product $ij$. Without loss of
generality, we may exhaust all possibilities by checking all of the
generators and zero.

If $ij=\pm 1$, then $iij=\pm i$, which implies that $j=\mp i.$
If $j$ can be written as a multiple of $i$, then they are not
linearly dependent. A basis for $\mathbb{D}$ is then the set $\{1,i\}$,
which contradicts our original assumption that the basis is $\{1,i,j\}$.

If $ij=\pm i$, then $iij=\pm i^{2}=\mp 1,$ which implies that $j=\pm 1.$
Again this would imply that $\mathbb{D}$ is a two-dimensional field extension over
$\mathbb{R}.$

Analogously, if $ij=\pm j,$then $ijj=\pm j^{2}=\mp 1,$which implies that $i=\pm 1,$
and we conclude as above.

Suppose finally that $ij=0.$ In this case there would exist a zero-divisor (in fact, both
$i$ and $j$ would be zero-divisors), which is incompatible with $\mathbb{H}$ being a field.
More directly, if $ij=0$, then we would have that $ijj=0,$ i.e. this would show that $i$
itself is zero, against the original assumption.

We have exhausted all possibilities and can conclude that there is no third
degree field extension over $\mathbb{R}$ that holds the property that $
i^{2}=j^{2}=-1$. Thus it is not possible to create the Theory of Triplets
while satisfying the requirements of a field.

\begin{remark} While our simple computations show the difficulty of building a three-dimensional field
extension over $\mathbb{R}$, this fact is in fact a special case of a much more powerful and general
result known as Frobenius theorem, \cite{palais}. In 1877, Ferdinand Georg Frobenius characterized the finite dimensional associative division algebras over the real numbers, and in particular proved that the only associative division algebra which is not commutative over the real numbers is the algebra of quaternions. In addition, Frobenius proved that
there are only three finite dimensional division algebras over $\mathbb{R}$: $\mathbb{R}$ itself, $\mathbb{C}$, and finally $\mathbb{H}.$
\end{remark}

This quick analysis shows why Hamilton had to consider a four
dimensional field extension by adding a new element k that is linearly
independent of the generators $1,i,$ and $j$. Thus our question becomes whether it is
possible to find a four-dimensional field extension of $\mathbb{R}$ whose basis is $\{1,j,i,k\}$
and such that $i^{2}=j^{2}=k^{2}=-1$. Specifically, the question is to discover (and
this must have been the issue that Hamilton was struggling with) what are the relationships
between the three units.

Since we have excluded all the previous possibilities on the outcome of $ij$, we will assume
that $ij=k$ and we will investigate the consequences of such an assumption.
To begin with, we want to understand what should the product $ji$ be.
By an argument analogous to the one above, we can reason that $ji$
yields either $k$ or $-k$. If $ji=k$, then $k^{2}=(ij)^{2}=i^{2}j^{2}=(-1)(-1)=1$.
This gives again a contradiction, because the product of two non-zero elements, $1+k$
and $1-k$ would then give $(1+k)(1-k)=1-k^{2}=1-1=0$. Once again we have zero-divisors
and $\mathbb{H}$ is not a field.

\begin{remark}
Note that $ij=ji$ is incompatible with the request that $\mathbb{H}$ is a field. However,
it is quite possible to consider such a structure, which is of independent interest
and is known as the commutative ring of bicomplex numbers \cite{Bicomplex}.
\end{remark}

We (and Hamilton!) are therefore left with the situation in which $ji=-k$. This generates a
non-commutative field and it is exactly the structure that Hamilton called with the name of quaternions.
It was a rather
radical move at the time to study a non-commutative group and this proposal was met by
many with skepticism \cite{Kuipers}. Most mathematicians, at the time, could not grasp the
value of studying an algebraic structure that violated the law of
multiplicative commutativity. Modern day mathematics recognize this as a
very valuable practice. American mathematician Howard Eves described the
significance of Hamilton's new approach by saying that it 
\begin{quote}
opened the
floodgates of modern abstract algebra. By weakening or deleting various
postulates of common algebra, or by replacing one of more of the postulates
by others, which are consistent with the remaining postulates, an enormous
variety of systems can be studied. As some of these systems we have groupoids, quasigroups, loops, semi-groups, monoids, groups, rings, integral domains, lattices, division rings, Boolean rings, Boolean algebras, fileds, vector spaces, Jordan algebras, and Lie algebras, the last two being examples non-associative algebras. \cite{Kuipers}.
\end{quote}

\bigskip

\section{Holomorphicity}

In this section, I will recall the notion of holomorphicity in complex analysis and
I will explore how such a notion can be extended to the case of quaternionic functions
defined on the space of quaternions. This is a very natural next step, since the power
and importance of complex numbers cannot be exploited until a full theory of holomorphic
functions is developed. For my discussion of holomorphic functions of a complex variable I
will follow \cite{Berenstein}.

\begin{definition} A function $f:\mathbb{C}\rightarrow \mathbb{C}$ is said to be complex
differentiable at $z_{0}\in \mathbb{C}$ if the limit $$\lim\limits_{z\rightarrow z_{0}}\frac{f(z)-f(z_{0})}{z-z_{0}}$$
exists, when $z$ approaches $z_0$ in the complex plane (the limit, therefore, is independent of the curve along which
$z$ approaches $z_0.$. This limit will be called the first derivative of $f$ at $z_0$ and will be denoted by
$f'(z_0).$
\end{definition}

\begin{definition} A function $f$ is said to be holomorphic in an open set $U\subset \mathbb{C},$ if
$f'(z)$ exists for all $z\in U$.
\end{definition}

In real analysis, i.e. for real valued functions of real variables, the property of differentiability is not
a very strong property. For example, there are plenty of functions which admit a first derivative, even a
continuous first derivative, and yet do not have a second derivative. There are also examples of functions which
admit infinitely many derivatives, and yet those derivatives are insufficient to reconstruct the original functions
by means of its Taylor series. In other words, there are functions which are $C^{\infty}$ but not real analytic. For this reason, the next result is a very surprising statement, and one that clearly indicates that the study of complex functions
is truly a different topic from real analysis.

\begin{theorem} Let $f:U\subset \mathbb{C} \to \mathbb{C}$ be a function. The following three properties
are equivalent:

1. $f$ admits a first derivative, in the complex sense, in every point of U.

2. $f$ is can be represented, near every point $\alpha \in U$, by a power series $\sum_{n=0}^{\infty
}c_{n}(z-\alpha )^{n}.$

3. $f$ is a solution of the Cauchy-Riemann equation on U. That is $$\frac{\partial f%
}{\partial x}+i\frac{\partial f}{\partial y}=0.$$
\end{theorem}
\bigskip

The first natural question that comes to mind, is whether these three conditions can be
reformulated in the quaternionic context, and whether they are still equivalent in that setting.
In a very well known paper, \cite{Sudbery}, A. Sudbery investigates and answers these questions.

To begin with, we need a definition
of quaternionic differentiability. Because we are working in a skew field,
we are not guaranteed commutativity, thus we must define derivatives from
the left and right side.

\begin{definition} A function is said to be quaternion-differentiable
on the left at a point $q\in \mathbb{H}$ if the limit $$\frac{df}{dq}%
=\lim\limits_{h\rightarrow 0}\frac{f(q+h)-f(q)}{h}$$ exists, when $h$ converges to
zero along any direction in the quaternionic space.
\end{definition}

It turns out, \cite{Sudbery}, that this is not a very good definition, since
if this limit exists, then $f$ must necessarily be of the form $f(q)=a+qb$ for $a,b\in
\mathbb{H}$. Thus if a function is left differentiable in the quaternionic
sense, it must be a linear function.

Given that differentiability is not the appropriate way to generalize the notion of
holomorphicity, one can investigate power series for quaternionic functions.

\begin{definition} A function $f:\mathbb{H}\rightarrow \mathbb{H}$ is said to be a quaternionic
monomial
if $f(q)=\Pi _{i=1}^{r}a_{i}q$ where $r$ is a non-negative number and $
a_{i}\in \mathbb{H}$.
\end{definition}

The reason for such a definition is due to the
non-commutativity of $\mathbb{H}$, which does not allow the coefficients $a_i$ to move to
the same side of the powers of $q$. It is easy to show, however, that the quaternionic
functions that can be represented as a power series of quaternionic
monomials are the complex functions that are real analytic, and so once again we have a definition that does
not offer any new class of functions.

So, if quaternionic differentiability confines us to linear functions, and expressibility
in power series limits us to real analytic functions, the only hope for an interesting
theory of holomorphicity seems to be based on the possibility of generalizing the Cauchy-Riemann
equations.

It was the Swiss mathematician Fueter, who was able to develop the appropriate generalization
of the Cauchy-Riemann equations to the quaternionic case.
The Cauchy-Fueter
equations (as they are now called in honor of their inventor) were not developed until 1935, almost a century after Hamilton's
discovery of quaternions. Because quaternionic derivatives are defined on
both the right and left side, there are actually two differential operators, the left
Cauchy-Fueter operator, $\frac{\partial _{l}}{\partial \overline{q}}$, and the
right Cauchy-Fueter operator, $\frac{\partial _{r}}{\partial \overline{q}}$.
These operators are defined as follows:

$$\bigskip \frac{\partial _{l}f}{\partial \overline{q}}=\frac{\partial f}{\partial
q_{o}}+i\frac{\partial f}{\partial q_{1}}+j\frac{\partial f}{\partial q_{2}}+k\frac{%
\partial f}{\partial q_{3}}$$

$$\frac{\partial _{r}f}{\partial \overline{q}}=\frac{\partial f}{\partial q_{o}}+%
\frac{\partial f}{\partial q_{1}}i+\frac{\partial f}{\partial q_{2}}j+\frac{\partial f%
}{\partial q_{3}}k.$$

These operators are used to define the Cauchy-Fueter equations, whose solutions are
the appropriate generalization of holomorphic functions to the quaternionic setting.
In accordance with the terminology established by Fueter in his early work, and now commonly adopted, these
functions are now known as regular functions.

Specifically, we have the following definition.

\begin{definition} Let $f:\mathbb{H} \to \mathbb{H}$ be a quaternionic valued function of a quaternionic
variable. We say that $f$ is left regular if
$$\frac{\partial _{l}f}{\partial \overline{q}}=\frac{%
\partial f}{\partial q_{o}}+i\frac{\partial f}{\partial q_{1}}+j\frac{\partial f}{%
\partial q_{2}}+k\frac{\partial f}{\partial q_{3}}=0$$

\noindent
and we say that $f$ is right regular if

$$\frac{%
\partial _{r}f}{\partial \overline{q}}=\frac{\partial f}{\partial q_{o}}+\frac{%
\partial f}{\partial q_{1}}i+\frac{\partial f}{\partial q_{2}}j+\frac{\partial f}{%
\partial q_{3}}k=0.$$

\end{definition}

As it is shown in \cite{Sudbery}, this definition turns out to be very effective, and
the theory of regular functions is a fully formed theory, which in many essential ways resembles the theory of holomorphic functions of a complex variable. For example, it is possible to reconstruct (almost with no
important changes) the entire Cauchy theory that holds for holomorphic functions (Cauchy theorems,
Cauchy representation formulas, etc.) and therefore almost all of its consequences.

An interesting observation is that any quaternionic valued function $f=f_0+f_1i+f_2j+f_3k$
can be thought of as a vector $\vec{f}=(f_0,f_1,f_2,f_3)$ of
four real valued functions, and therefore a regular function corresponds to a vector $\vec{f}$
of real functions satisfying a $4\times 4$ system of linear constant coefficients partial differential equations, whose
formal appearance is the same we saw when we represented quaternionic multiplication via a matrix (we leave this simple computation to the reader).

There are, nevertheless, some important differences that emerge between the two theories, and some
are substantial enough to justify the emergence of a new approach. I want to conclude this
historical review by pointing out the nature of the differences, and one of the new approaches that
seem to remedy and offer a different theory.

To begin with, it is an unpleasant discovery that even $f(q)=q=q_0+q_1i+q_2j+q_3k$ is not regular according to
the definition of Fueter. Indeed,

$$\frac{\partial _{r}q}{\partial \overline{q}}=\frac{\partial q}{\partial q_{o}}+\frac{%
\partial q}{\partial q_{1}}i+\frac{\partial q}{\partial q_{2}}j+\frac{\partial q}{%
\partial q_{3}}k=1-1-1-1=-2\neq 0.$$

Regardless of what one has in mind with this new theory, the fact that the simplest function is
not regular, and therefore that no monomial or polynomial is regular is problematic.

There have been several attempts to circumvent this issue, but the most promising, and recent,
has just appeared in \cite{advances}, where the following definition is offered.

\begin{definition} Let $U$ be a domain in $\mathbb{H}.$ A real differentiable function $f:U \to \mathbb{H}$ is said to
be slice-regular if, for every quaternion $I=x_1i+x_2j+x_3k$ such that $x_1^2+x_2^2+x_3^2=1$, the
restriction of $f$ to the complex line $L_I=\mathbb{R}+\mathbb{R}I$ passing through the origin, and containing $1$ and $I$
is holomorphic on $U \cap L_I.$
\end{definition}

Note that this definition essentially requires that a function be holomorphic on each complex slice of the original
domain: this explain the terminology of 'slice-regularity.' In a sense, this does not seem to be such an
interesting idea, but the consequences of this definition are quite far reaching. For example, one can prove that the
identity function $f(q)=q$ is slice-regular but, more important, all
monomials of the form $aq^n$, with $a$ any quaternion, are slice-regular in this sense. Since regularity respects addition, it is immediate to see that all polynomials of the form $f(q)=\Sigma_{t=0}^n a_nq^n$ with $a_n \in \mathbb{H}$ are regular and
in fact even power series of the form $f(q)=\Sigma_{t=0}^{\infty} a_nq^n$, with $a_n \in \mathbb{H}$ are regular where convergent, \cite{advances}.

This theory is only in its infancy at this point (though a cursory check to mathsci.net shows that a significant number of
papers have already appeared in this area), but it seems to demonstrate that the history of quaternions, which began in 1835, is not yet concluded, and that new ideas are constantly germinating in this field.

\bigskip \textbf{\pagebreak }

\section{Appendix: Letter from Sir W.R.Hamilton to Rev. Archibald H. Hamilton}

Letter dated August 5, 1865.

MY DEAR ARCHIBALD - (1) I had been wishing for an occasion of corresponding
a little with you on QUATERNIONS: and such now presents itself, by your
mentioning in your note of yesterday, received this morning, that you
\textquotedblleft have been reflecting on several points connected with
them\textquotedblright\ (the quaternions), \textquotedblleft particularly on
the Multiplication of Vectors.\textquotedblright\ (2) No more important, or
indeed fundamental question, in the whole Theory of Quaternions, can be
proposed than that which thus inquires What is such MULTIPLICATION? What are
its Rules, its Objects, its Results? What Analogies exist between it and
other Operations, which have received the same general Name? And finally,
what is (if any) its Utility? (3) If I may be allowed to speak of myself in
connexion with the subject, I might do so in a way which would bring you in,
by referring to an ante-quaternionic time, when you were a mere child, but
had caught from me the conception of a Vector, as represented by a Triplet:
and indeed I happen to be able to put the finger of memory upon the year and
month - October, 1843 - when having recently returned from visits to Cork
and Parsonstown, connected with a meeting of the British Association, the
desire to discover the laws of the multiplication referred to regained with
me a certain strength and earnestness, which had for years been dormant, but
was then on the point of being gratified, and was occasionally talked of
with you. Every morning in the early part of the above-cited month, on my
coming down to breakfast, your (then) little brother William Edwin, and
yourself, used to ask me, \textquotedblleft Well, Papa, can you multiply
triplets\textquotedblright ? Whereto I was always obliged to reply, with a
sad shake of the head: \textquotedblleft No, I can only add and subtract
them.\textquotedblright\ (4) But on the 16th day of the same month - which
happened to be a Monday, and a Council day of the Royal Irish Academy - I
was walking in to attend and preside, and your mother was walking with me,
along the Royal Canal, to which she had perhaps driven; and although she
talked with me now and then, yet an under-current of thought was going on in
my mind, which gave at last a result, whereof it is not too much to say that
I felt at once the importance. An electric circuit seemed to close; and a
spark flashed forth, the herald (as I foresaw, immediately) of many long
years to come of definitely directed thought and work, by myself if spared,
and at all events on the part of others, if I should even be allowed to live
long enough distinctly to communicate the discovery. Nor could I resist the
impulse - unphilosophical as it may have been - to cut with a knife on a
stone of Brougham Bridge, as we passed it, the fundamental formula with the
symbols, i, j, k; namely, $i^{2}=j^{2}=k^{2}=ijk=-1$which contains the
Solution of the Problem, but of course, as an inscription, has long since
mouldered away. A more durable notice remains, however, on the Council Books
of the Academy for that day (October 16th, 1843), which records the fact,
that I then asked for and obtained leave to read a Paper on Quaternions, at
the First General Meeting of the session: which reading took place
accordingly, on Monday the 13th of the November following.

With this quaternion of paragraphs I close this letter I.; but I hope to
follow it up very shortly with another.

Your affectionate father,

WILLIAM ROWAN HAMILTON.

\bigskip

\pagebreak

\end{document}